%% file: main_tex.tex
\documentclass[preprint, twocolumn,5p]{elsarticle}

\usepackage{amssymb,amsmath,amsfonts}
\usepackage{latexsym}
\usepackage{graphicx}
\usepackage{color}
\usepackage{tikz}

\newcommand{\N}{\mathbb{N}}
\newcommand{\R}{\mathbb{R}}
\newcommand{\PP}{\mathbb{P}}
\newcommand{\EE}{\mathbb{E}}
\newcommand{\e}{\varepsilon}
\newcommand{\dd}{{\rm d}}
\newcommand{\oo}{{\rm o}}

\newtheorem{thm}{Theorem}
\newtheorem{prop}{Proposition}
\newtheorem{assumption}{Assumption}
\newtheorem{lem}[thm]{Lemma}
\newdefinition{rmk}{Remark}
\newproof{pf}{Proof}
\newproof{pot}{Proof of Theorem \ref{thm2}}

\journal{~}

\begin{document}

\begin{frontmatter}

\title{A multi-scale study of a class of hybrid predator-prey models}

\author{Alexandre Genadot}

\address{Centre De Recherche en Math\'ematiques de la D\'ecision,\\ Universit\'e PARIS-DAUPHINE,\\ Place du Maréchal De Lattre De Tassigny,\\ 75775~PARIS~CEDEX~16, FRANCE }

\ead{Alexandre.Genadot@math.cnrs.fr}

\begin{abstract}
We address the question of an averaging principle for a general class of multi-scale hybrid predator-prey models. We consider prey-predator models where the kinetic of the prey population, described by a differential equation, is faster than the kinetic of the predator population, described by a jump process, the two dynamics being fully coupled. An averaged model is obtained and compared to the original slow-fast system in term of probability and absorption time.
\end{abstract}


\end{frontmatter}

\section{Introduction}

We study a general class of multi-scale hybrid models of predator-prey type. These models describe the interactions between two populations of individuals at a mesoscopic scale: the prey population is assumed to be large enough so that a deterministic approximation can be justified, its dynamic following a deterministic differential equation, whereas the predator population evolves according to a stochastic jump process, see (\ref{intro:eq1}). An other way to describe such a hybrid situation is to consider that the prey population evolves according to a fast time scale so that its evolution is somehow averaged and thus described by a continuous deterministic equation whereas the dynamic of the predator population is slower so that stochastic features remain and may be therefore best described by a stochastic jump process.
\begin{equation}\label{intro:eq1}
\left\{\begin{array}{ll}
{(x_t)}_{t\geq0}& \text{preys: continuous deterministic dynamic,}\\
{(n_t)}_{t\geq0}& \text{predators: discrete stochastic dynamic.}
\end{array}\right.
\end{equation}

Such hybrid models are widely used in biology and may describe a large range of situations. For instance, in neurosciences, the generation of an action potential is modeled by a differential equation (or a partial differential equation) fully coupled to a jump mechanism accounting for the conformal changes of the ion-channels, see \cite{BR11}. In a resource-consumption situation, hybrid models have been used to describe the evolution of a bacteria population in bio-reactors, see \cite{CO79,CMM13}. As explained in \cite{SGM07} with examples coming from cancerology, hybrid models allow to describe biological phenomena which are in interactions but evolving on different space or time scales.

Indeed, a common feature of many biological systems is to be intrinsically endowed with different time scales. We place ourself in this framework considering that the prey dynamic is faster than the predator one. The resulting model is then a slow-fast hybrid model that we intend to reduce through the averaging principle. Reduction of slow-fast hybrid models is a quite recent field of research, the first mathematical analysis being, as far as we know, the work \cite{Fa10}, followed from various works of the authors of \cite{Wa12b}. Get reduced models from slow-fast ones allows to simplify the equations describing the system of interest and in this way, this allows to perform a more detailed analysis of the biological phenomena in question, both theoretically and numerically. For instance, in our case, the probability and time of extinction for the predator population are more accessible on the averaged model than on the slow-fast one.

Compared to the above mentioned studies \cite{Fa10,Wa12b} or to the works \cite{Ge12,Wa12} on multi-scale hybrid systems, the present work has the originality to propose to average a slow-fast hybrid model with respect to its continuous component. In the previously mentioned works, it was natural to consider averaging with respect to the discrete component of the system. We believe that the case that we consider in the present paper is as much relevant in some biological framework where the continuous deterministic variable may evolve on a faster time scale than the  discrete stochastic variable, see for instance Section \ref{sec:ts}. Let us also remark, even if it is of secondary importance, that we work in the present paper with discrete components which are countable and not only finite contrary to the aforementioned studies.

The paper is organized as follows. In Section \ref{bsec:mod}, we describe the model that we are interested in. We motivate our presentation by the description of a very particular, but quite relevant, example in Section \ref{sec:HPP}. We endow this example with different time scales in Section \ref{sec:ts}. The general class of slow-fast hybrid models considered in the present paper is then described in Section \ref{sec:gs}. In Section \ref{bsec:main}, our main result, the averaging principle for this class of processes, is introduced together with some important properties of the averaged models. We present numerical simulations to illustrate the obtained results in Section \ref{bsec:ex}. The \ref{sec:pf} is devoted to the proof of the main result of Section \ref{bsec:main}.

\section{The model}\label{bsec:mod}

\subsection{A hybrid predator-prey models}\label{sec:HPP}

\begin{figure}
\begin{center}
\begin{tikzpicture}
\node[draw,circle,fill=gray!50] (A) at (0,0) {$0$};
\node[draw,circle] (B) at (1.5,0) {$1$};
\node[draw,circle] (C) at (3,0) {$2$};
\node[draw,circle] (D) at (4.5,0) {$3$};
\node[draw,circle] (E) at (6,0) {$4$};
\node[draw,circle,dashed] (F) at (7.5,0) {\textcolor{white}{$5$}};
\draw[->,>=latex] (B) to[bend left=45] node[midway, above]{$\mu(\zeta)$} (C);
\draw[->,>=latex] (C) to[bend left=45] node[midway, above]{$2\mu(\zeta)$} (D);
\draw[->,>=latex] (D) to[bend left=45] node[midway, above]{$3\mu(\zeta)$} (E);
\draw[->,>=latex,dashed] (E) to[bend left=45] (F);
\draw[->,>=latex] (B) to[bend right=-45] node[midway, below]{$D$} (A);
\draw[->,>=latex] (C) to[bend right=-45] node[midway, below]{$2D$} (B);
\draw[->,>=latex] (D) to[bend right=-45] node[midway, below]{$3D$} (C);
\draw[->,>=latex] (E) to[bend right=-45] node[midway, below]{$4D$} (D);
\draw[->,>=latex,dashed] (F) to[bend right=-45]  (E);
\end{tikzpicture}
\caption{Kinetic of the jump process ${(n_t)}_{t\geq0}\subset \N_0$ of Setion \ref{sec:HPP} if ${(x_t)}_{t\geq0}$ were held fixed to some real $\zeta$ and $\beta=\gamma=1$. In gray, $0$ is the absorbing state.}\label{fig:kin:jump}
\end{center}
\end{figure}
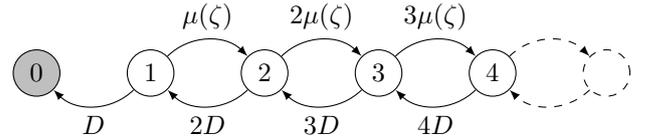

Let us consider a population of preys and predators in some domain $V$ (area or volume). The number of preys in the population at time $t$ is denoted by $x_t$ while $y_t$ denotes the number of predators. Suppose, in a very simple model, that the predators and preys die respectively at rate $\gamma D$ and $D$, $\gamma$ being some positive ratio parameter. The growth of the predator population is assumed to depend on the number of preys through the rate $\beta\mu(\cdot)$, $\beta$ being some positive conversion efficiency and $\mu$ some consumption (of preys) function. Accordingly, the prey population decreases proportionally to the predator population at rate $\frac{\alpha}{V}\mu(\cdot)$, where $\alpha$ may also be seen as a conversion efficiency (or as the inverse of the yield coefficient in a bacteria population). In such a situation, to give a chance to the prey population to survive, one may add some immigration of preys at constant flow $Dx_{in}$ with $x_{in}$ a positive parameter. This model may be described by the two following differential equations
\begin{equation}\label{eq:HPP:det}
\left\{\begin{array}{rl}
\frac{\dd x_t}{\dd t}&=D(x_{in}-x_t)-\frac{\alpha}{V}\mu(x_t)y_t,\\
\frac{\dd y_t}{\dd t}&=(\beta\mu(x_t)-\gamma D)y_t,
\end{array}\quad t\geq0,\right.
\end{equation}
endowed with some initial condition $(x_0,y_0)$. In this model, only the predator population responds to primary production of preys. This kind of prey-dependent models are particularly appropriate for so-called homogeneous systems like bacteria feeding in a stirred chemostat, see \cite{CO79,CMM13,AG89}, which may also be seen as a situation of resource-consumption (the preys becoming the resource and the predators the consumers). More general situations are included in the more general setting presented in Section \ref{sec:gs}.\\
We actually made a choice when describing the model by the system of differential equations (\ref{eq:HPP:det}). We choose to describe the model at a macroscopic level at which we don't have to describe the possible discrete jump of population individuals. We implied for example that the prey population was large enough so that a deterministic approximation can be justified. At a mesoscopic scale, if the predator population at time $t$ is composed of $n$ individuals, it looses one of them at rate $n\gamma D$ (death of one of the $n$ predators) or on the contrary gains one individual at rate $n\beta \mu(x_t)$ (one of the $n$ predators give birth). The model is then more appropriately described by the equation
\begin{equation}\label{eq:x}
\frac{\dd x_t}{\dd t}=D(x_{in}-x_t)-\frac{\alpha}{V}\mu(x_t)n_t,\quad t\geq0,
\end{equation}
fully coupled to the jump mechanism
\begin{align}
&\PP(n_{t+h}=n_2|n_t=n_1)\nonumber\\
=&\left\{\begin{array}{ll}
\beta \mu(x_t)n_1h+\oo(h)&\textrm{if } n_2=n_1+1,n_1\geq1,\\
\gamma Dn_1h+\oo(h)&\textrm{if } n_2=n_1-1,n_1\geq1,\\
1-(\beta \mu(x_t)+\gamma D)n_1\hspace{-0.3cm}&h +\oo(h)\\
&\textrm{if } n_2=n_1,n_1\geq1,\\
\oo(h)&\text{else},
\end{array}\label{eq:n}\right.
\end{align}
endowed with some initial condition $(x_0,n_0)$ (possibly random). The integer-valued\footnote{The notation $\N_0$ stands for the set of non-negative integers.} process ${(n_t)}_{t\geq0}\subset \N_0$ describes the evolution of the predator population at a discrete level. This is a jump process absorbed in zero. If ${(x_t)}_{t\geq0}$ were held fixed to some real $\zeta$, the kinetic of ${(n_t)}_{t\geq0}$ would be the one of a homogeneous birth and death process absorbed in zero with parameters $\beta\mu(\zeta)$ and $\gamma D$ \cite[see e.g. on birth and death processes][Ch. 6]{PP10}. The kinetic of this jump process is illustrated in Figure \ref{fig:kin:jump}.

\noindent The process described by equation (\ref{eq:x}-\ref{eq:n}) is a so-called hybrid process combining a macroscopic continuous deterministic behavior fully coupled to a microscopic discrete stochastic behavior. Mathematically, assuming that $\mu$ is continuous, one may show that the couple ${(x_t,n_t)}_{t\geq0}$ is a piecewise deterministic Markov process in the sens of \cite{D93}. Its generator $\mathcal{A}$ describes the infinitesimal motion of the process
\begin{align*}
\mathcal{A}\phi(x,n)=&~\left[D(x_{in}-x)-\frac{\alpha}{V}\mu(x)n\right]\phi'(x,n)\\
&+[\phi(x,n+1)-\phi(x,n)]\beta \mu(x)n\\
&+[\phi(x,n-1)-\phi(x,n)]\gamma Dn,
\end{align*}
defined for any bounded function $\phi:\R\times\N_0\to\R$ continuously derivable in its first variable, measurable in its second variable. Notice that the first part of the generator describes the deterministic motion of the process between two jumps, the jumps being described by the second part of the generator.

\begin{figure}
\begin{center}
\scalebox{.6}{\input{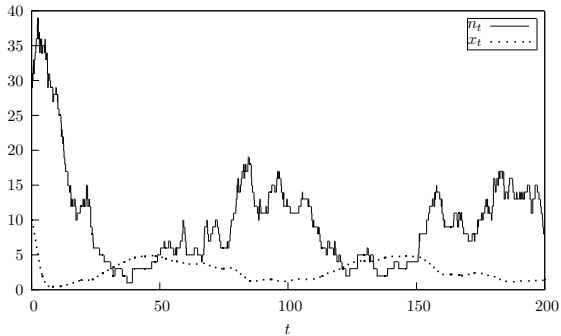}}
\caption{We present a trajectory of the process ${(x_t,n_t)}_{t\in[0,200]}$ of Section \ref{sec:HPP} satisfying equations (\ref{eq:x}-\ref{eq:n}) with initial condition $(x_0,n_0)=(10,30)$. The dotted curve is for the process ${(x_t)}_{t\in[0,200]}$ while the piecewise constant curve is for ${(n_t)}_{t\in[0,200]}$, both curves being drawing with respect to time. We choose the Monod function $\mu(x)=\frac{0.15x}{1+x}$ as consumption function. The other parameter values are: $D=0.1$, $V=1$, $x_{in}=7$, $\alpha=0.5$, $\beta=\gamma=1$.}\label{fig:traj}
\end{center}
\end{figure}

\noindent Figure \ref{fig:traj} displays a trajectory of the process ${(x_t,n_t)}_{t\geq0}$. For efficient numerical simulations of piecewise deterministic process, we refer to \cite{R13}.

\noindent It is quite natural to question the link between the two models (\ref{eq:HPP:det}) and (\ref{eq:x}-\ref{eq:n}). The answer is that the first one may be obtained as a scaling limit of the second one. This scaling limit has to be understood in the following way: let us renormalized the volume and the number of predators by a factor $N\in\N$ in equation (\ref{eq:x}-\ref{eq:n}): $V\leadsto NV$, $n\leadsto\frac{1}{N}n$. The renormalized process ${(x^N_t,n^N_t)}_{t\geq0}$ still satisfies equations (\ref{eq:x}-\ref{eq:n}) except that the process ${(n^N_t)}_{t\geq0}$ is valued in $\frac{1}{N}\N_0$ and jumps at rates multiply by $N$. The following result is derived from \cite{CJM14}.
\begin{thm}\cite[in][Theorem 2.3]{CJM14}
For any time horizon $T$, the rescaled hybrid process ${(x^N_t,n^N_t)}_{t\in[0,T]}$ converges in law in\footnote{$\mathcal{D}([0,T],\R_+)$ is the Skorohod space of real-valued c\`adl\`ag function on $[0,T]$.} $\mathcal{C}([0,T],\R_+)\times\mathcal{D}([0,T],\R_+)$ when $N$ goes to infinity, towards the deterministic process ${(x_t,y_t)}_{t\in[0,T]}$ solution of (\ref{eq:HPP:det}).
\end{thm}
This theorem provides a way to reduce the complexity of the hybrid model: in the framework of large volume, one may use the deterministic model solution of (\ref{eq:HPP:det}) as an approximation of the hybrid model (\ref{eq:x}-\ref{eq:n}). One may argue at this point that with such an approximation, one looses the intrinsic variability of the hybrid model and potential "finite size" effects that may arise in discrete models. To recover some variability, one may study the fluctuation of the rescaled hybrid model around its deterministic limit, that is look for a central limit theorem. Another approach consists in taking advantage of the different time scales -- if there is indeed such different time scales -- of the hybrid model components to simplify its dynamic. The present paper develop this second approach.

\subsection{A two time-scale model.}\label{sec:ts}

In this section, we introduce a time scale in the model (\ref{eq:x}-\ref{eq:n}). Suppose that both the conversion efficiency coefficient $\beta$ and the death rate coefficient $\gamma$ for the predators are low. Mathematically, this consideration leads to introduce a (small) parameter $\e$ in the model, replacing $\beta$ by $\e\beta$ and $\gamma$ by $\e\gamma$ in equation (\ref{eq:n}). The quantity $\e$ is a time scale parameter: the introduction of $\e$ has the effect of slowing down the rate at which the process jumps. Thus, the jumping part (predators) of the process evolves on a slower time scale than the continuous component (preys) of the process. Conversely, on can say that the continuous component evolves on a faster time scale than the jumping one. This is this formulation which leads to write the slow-fast model in an appropriate form for slow-fast analysis. We denote by ${(x^\e_t,n^\e_t)}_{t\geq0}$ the two-time scale process. It satisfies the differential equation
\begin{equation}\label{eq:xe}
\frac{\dd x^\e_t}{\dd t}=\frac{1}{\e}\left[D(x_{in}-x^\e_t)-\frac{\alpha}{V}\mu(x^\e_t)n^\e_t\right],\quad t\geq0,
\end{equation}
fully coupled to the jump mechanism
\begin{align}
&\PP(n^\e_{t+h}=n_2|n^\e_t=n_1)\nonumber\\
=&\left\{\begin{array}{ll}
\beta \mu(x^\e_t)n_1h+\oo(h)&\textrm{if } n_2=n_1+1,n_1\geq1,\\
\gamma Dn_1h+\oo(h)&\textrm{if } n_2=n_1-1,n_1\geq1,\\
1-(\beta \mu(x^\e_t)+\gamma D)n_1\hspace{-0.3cm}&h +\oo(h)\\
&\textrm{if } n_2=n_1,n_1\geq1,\\
\oo(h)&\text{else},
\end{array}\label{eq:ne}\right.
\end{align}
with initial condition $(x^\e_0,n^\e_0)$ (possibly random). Slow-fast systems are very common in biology. Besides the considered population dynamic framework, they arises for instance naturally in neurosciences \cite{Wa12} or molecular biology \cite{Fa10}. A good starting point to learn multi-scale techniques certainly is the excellent book \cite{PS08}. For two-time scale methods applied to slow-fast continuous Markov chains, we refer to \cite{YZ12}.\\
Our aim is to reduce the complexity of the model (\ref{eq:xe}-\ref{eq:ne}) by taking advantage of the two different time scales. The reduced model is expected to be easier to handle both theoretically and numerically. This provides a way to go further in the understanding of the distribution and the structure of the underlying uncertainty of the hybrid model (\ref{eq:xe}-\ref{eq:ne}).\\
The reduction we intend to perform relies on averaging techniques: the dynamic of the component $x^\e$ of the process is accelerated due to the $\frac1\e$ scaling. Letting formally $\e$ goes to zero, we accelerate the dynamic of $x^\e$ so much that it will reach instantaneously its stationary behavior (if it has one). Thus, the dynamic of the slow component $n^\e$ will be somehow averaged against the stationary behavior of the fast component.  The next sections provide a rigorous framework to the above heuristic.

\subsection{General setting}\label{sec:gs}

The model presented in the previous section belongs to a more general class of hybrid models that we now proceed to define. For a time scale parameter $\e\in(0,1)$, let ${(x^\e_t,n^\e_t)}_{t\geq0}$ be the hybrid process defined on some probability space $(\Omega,\mathcal{F},\PP)$ with generator $\mathcal{A^\e}$ given by
\begin{align}\label{g:main}
\mathcal{A^\e}\phi(x,n)=&~\frac1\e g(x,n)\phi'(x,n)\nonumber\\
&+[\phi(x,n+1)-\phi(x,n)]b(x,n)\\
&+[\phi(x,n-1)-\phi(x,n)]d(x,n),\nonumber
\end{align}
defined for any bounded function $\phi:\R\times\N_0\to\R$ continuously derivable in its first variable, measurable in its second variable. The following assumption gathers some basic hypotheses about the main characteristic of the model, the functions $g$, $b$ and $d$. 
\begin{assumption}\label{ass:2}
The functions $g:\R\times\N_0\to\R$ and $b,d:\R\times\N_0\to\R_+$ are assumed to be continuous in their first variable and measurable in their second variable. Moreover, $b(\cdot,0)=d(\cdot,0)=0$ but otherwise $b$ and $d$ are positive on $\R\times\N$.
\end{assumption}
In a more dynamic view, the form of the generator (\ref{g:main}) together with Assumption \ref{ass:2} mean that the process ${(x^\e_t,n^\e_t)}_{t\geq0}$ satisfies the differential equation
\begin{equation}\label{eq:gx}
\frac{\dd x^\e_t}{\dd t}=\frac{1}{\e}g(x^\e_t,n^\e_t),\quad t\geq0,
\end{equation}
fully coupled to the jump mechanism
\begin{align}
&\PP(n^\e_{t+h}=n_2|n^\e_t=n_1)\nonumber\\
=&\left\{\begin{array}{ll}
b(x^\e_t,n_1)h+\oo(h)&\textrm{if } n_2=n_1+1,n_1\geq1,\\
d(x^\e_t,n_1)h+\oo(h)&\textrm{if } n_2=n_1-1,n_1\geq1,\\
1-(b(x^\e_t,n_1)+d(x^\e_t,n_1))\hspace{-0.3cm}&h +\oo(h)\\
&\textrm{if } n_2=n_1,n_1\geq1,\\
\oo(h)&\text{else},
\end{array}\label{eq:gn}\right.
\end{align}
endowed with some initial condition $(x^\e_0,n^\e_0)\in\R_+\times\N_0$ (possibly random). Notice again that if ${(x^\e_t)}_{t\geq0}$ were held fixed to some real $\zeta$, the kinetic of ${(n^\e_t)}_{t\geq0}$ would be the one of a homogeneous birth and death process absorbed at zero with birth and death rates $b(\zeta,n)$ and $d(\zeta,n)$ ($n\in\N$) respectively with the absorbing condition $b(\zeta,0)=0$.
\begin{rmk}
The model (\ref{eq:gx}-\ref{eq:gn}) naturally includes the model described in the previous section. However, it is much more general and includes in particular so-called ratio-dependent models, such as those described in \cite{AG89}. In such models, the trophic function $\mu$ is not only determined by the prey-abundance ($\mu(x)$) but is rather a function of the prey abundance \textit{per capita} ($\mu\left(\frac xn\right)$).
\end{rmk}
\begin{assumption}\label{hyp:growth}
We assume the following growth conditions on the rate function $b$ and $d$. For any positive real $K$, there exist $c_1\geq0$ and $c_2,c_3>0$ three constants such that for any $n\in\N$
$$
\sup_{x\in[0,K]} b(x,n)+d(x,n)\leq c_1+c_2n,
$$
and for any $x,y\in\R$
\begin{align*}
&|b(x,n)-b(y,n)|+|d(x,n)-d(y,n)|\\
&\leq c_1(1+c_2n+c_3n^2)|x-y|.
\end{align*}
\end{assumption}
This assumption means that the rate functions $b$ and $d$ satisfy a homogeneous in $x$ sub-linear growth condition in $n$ and a non-homogeneous in $n$ Lipschitz condition in $x$. For the example developed in Section \ref{sec:HPP}, we have, for $(x,n)\in\R\times\N_0$, $b(x,n)=\beta\mu(x)n$ and $d(x,n)=Dn$. As long as $\mu$ is bounded and Lipschitz (which is the case for Monod functions for example, see the caption of Figure \ref{fig:traj}), Assumption \ref{ass:2} is satisfied.
\begin{assumption}\label{hyp:1}
We assume the three following conditions on the function $g$:
\begin{enumerate}
\item The function $g$ is strictly dissipative with rate $\delta$ in its first variable, uniformly in $n\in\N_0$. There exists $\delta>0$ such that for all $x_1,x_2\in \R$ and $n\in\N_0$
\begin{equation}\label{eq:diss}
(x_1-x_2)(g(x_1,n)-g(x_2,n))\leq -\delta(x_1-x_2)^2.
\end{equation}
\item The sequence ${(g(0,n))}_{n\in\N_0}$ is positive and bounded. 
\item For any $n\in\N_0$, the equation in $x$
$$
g(x,n)=0
$$
has a unique positive solution denoted $x^*_n$ depending at least measurably on $n$.
\end{enumerate}
\end{assumption}
Since the variable $x$ evolves in dimension $1$, the dissipativity condition \ref{hyp:1}.1 simply means that the function $g(\cdot,n)$ is strictly decreasing for all $n\in\N_0$. From Assumption \ref{hyp:1}, one can also show that the sequence of stationary points ${(x^*_n)}_{n\in\N_0}$ is bounded. Indeed, for $n\in\N_0$, taking $x_1=x^*_n$ and $x_2=0$ in (\ref{eq:diss}), we obtain $x^*_n\leq\frac{1}{\delta}g(0,n)\leq \frac{1}{\delta}\sup_{n\in\N_0}g(0,n)$. The fact that the sequence ${(g(0,n))}_{n\in\N_0}$ is positive also implies that if $x^\e_0>0$, then $x^\e_t>0$ for all $t\geq0$, $\PP$-a.s. Let us emphasize these facts in a proposition. 
\begin{prop}\label{prop:bound}
The two following statements hold,
\begin{itemize}
\item The hybrid process defined by the equations (\ref{eq:gx}-\ref{eq:gn}) is well defined. Moreover, if $x^\e_0>0$ $\PP$-a.s., then $x^\e_t>0$ for all $t\geq0$, $\PP$-a.s.
\item The sequence of stationary points ${(x^*_n)}_{n\in\N_0}$ is bounded,
$$
\sup_{n\in\N_0}x^*_n\leq \frac{1}{\delta}\sup_{n\in\N_0}g(0,n).
$$
\end{itemize}
\end{prop}
Notice that, as long as $\mu$ is a bounded non-negative strictly dissipative function such that $\mu(0)=0$, Assumption \ref{hyp:1} is satisfied for the example in Section \ref{sec:HPP} where for $(x,n)\in\R\times\N_0$, $g(x,n)=D(x_{in}-x)-\frac{\alpha}{V}\mu(x)n$.

\section{Main results}\label{bsec:main}

In this section, our aim is to reduce the complexity of the system described by the equations (\ref{eq:gx}-\ref{eq:gn}) by taking advantage of the presence of two time scales. As already explained at the end of section \ref{sec:ts}, heuristically, the picture is as follows. Letting $\e$ go to zero, we accelerate the dynamic of the fast variable $x^\e$ so that it instantaneously reaches its stationary behavior. Then, we average the dynamic of the slow variable $n^\e$ with respect to this stationary behavior. The resulting averaged process is described by one equation only and thus, is also called the reduced model. The dynamic of the reduced model is expected to be easier to handle than the dynamic of the slow-fast system, both theoretically and numerically. We now proceed to the rigorous statement of our averaging result.\\ 
Since for any fixed $n\in\N_0$ the fast subsystem $\frac{\dd x_t}{\dd t}=\frac{1}{\e}g(x_t,n)$, $(t\geq0)$ is dissipative according to Assumption \ref{hyp:1}, it will consume some energy until it reach its quasi-equilibrium $x^*_n$. The term "quasi" stands here to emphasize the fact that the equilibrium $x^*_n$ actually depends on the extra variable $n$. The proposition below precises at which rate this quasi-equilibrium is reached.
\begin{prop}\label{cv:fast}
For any $n\in\N_0$, consider the fast subsystem satisfying the differential equation
$$
\frac{\dd x_t}{\dd t}=\frac{1}{\e}g(x_t,n),\quad t\geq0,
$$
endowed with some positive initial condition $x_0$. Then the process ${(x_t)}_{t\geq0}$ satisfy
\begin{equation}
|x_t-x^*_n|\leq |x_0-x^*_n| e^{-\frac t\e\delta}.
\end{equation}
\end{prop}
\begin{pf}
Recall that, according to Assumption \ref{hyp:1}, the function $g$ is strictly dissipative with rate $\delta$ and, by definition of $x^*_n$, $g(x^*_n,n)=0$ for any $n\in\N_0$. Using these two facts, one may write
\begin{align*}
\frac{\dd}{\dd t}|x_t-x^*_n|^2&=\frac2\e(x_t-x^*_n)g(x_t,n)\\
&=\frac2\e(x_t-x^*_n)(g(x_t,n)-g(x^*_n,n))\\
&\leq -\frac2\e\delta|x_t-x^*_n|^2.
\end{align*}
This yields the result.
\end{pf}
Since, according to Proposition \ref{prop:bound}, $\sup_{n\in\N_0}|x^*_n|<\infty$, the latter proposition means that the isolated fast subsystem converges exponentially fast towards its equilibrium uniformly in $n\in\N_0$. This uniform convergence is crucial to prove the following averaging result.
\begin{thm}\label{thm:main}
Assume that, uniformly in $\e$, both $\EE((n^\e_0)^2)$ and the support of the law of $x^\e_0$ -- included in $\R_+$ -- are bounded. For any finite time horizon $T$, the jump process ${(n^\e_t)}_{t\in[0,T]}$ converges in law in $\mathcal{D}([0,T], \N_0)$ when $\e$ goes to zero, towards the averaged jump process ${(\bar{n}_t)}_{t\in[0,T]}$ with generator
$$
\mathcal{A}\phi(n)=[\phi(n+1)-\phi(n)]b(x^*_n,n)+[\phi(n-1)-\phi(n)]d(x^*_n,n),
$$
defined for any bounded measurable real function $\phi$ defined on $\N_0$.
\end{thm}
\begin{pf}
The proof is postponed to \ref{sec:pf}.
\end{pf}
\begin{rmk}
It is worth noticing that one may deduce from Theorem \ref{thm:main} -- and under the same assumptions -- that, for any finite time horizon $T$, the process ${(x^\e_t)}_{t\in[0,T]}$ converges in law in $\mathcal{D}([0,T], \R_+)$ when $\e$ goes to zero, towards ${(x^*_{\bar{n}_t})}_{t\in[0,T]}$.
\end{rmk}
According to this theorem, the averaged process is a homogeneous in time birth and death process on $\N_0$ with respective parameters
$$
\bar b_n=b(x^*_n,n),\quad\textrm{and}\quad \bar d_n=d(x^*_n,n),
$$
for $n\in\N$ and $b_0=0$ so that $0$ is still an absorbing state. The dynamic of the averaged process is easier to handle than the one of the slow-fast model. Indeed, birth and death processes with $0$ as absorbing state have been extensively studied in the literature \cite[see e.g.][Ch. 6]{PP10}. Let us mention, for example, some results about the probability and time of absorption for such a process.\\
Let $p_m$ and $t_m$ ($m\in\N$) denote respectively the probability of absorption into state $0$ and the mean absorption time, starting from some initial state $m$ in both cases.
\begin{thm}\cite[in][Theorem 6.1]{PP10}\label{thm:abs}
The absorption probability into state $0$ from the initial state $m$ is
\begin{equation}
p_m=\left\{\begin{array}{cl}
\frac{\sum_{i=m}^\infty\rho_i}{1+\sum_{i=1}^\infty\rho_i}&\textrm{if}\quad \sum_{i=1}^\infty\rho_i<\infty,\\
1&\textrm{if}\quad \sum_{i=1}^\infty\rho_i=\infty.
\end{array}\right.
\end{equation}
The mean absorption time is
\begin{equation}
t_m=\left\{\begin{array}{cl}
\sum_{i=1}^\infty\frac1{\bar b_i\rho_i}+&\hspace{-0.3cm}\sum_{k=1}^{m-1}\rho_k\sum_{j=k+1}^{\infty}\frac{1}{\bar b_j\rho_j},\\
~&\textrm{if}\quad \sum_{i=1}^\infty\frac1{\bar b_i\rho_i}<\infty,\\
\infty,&\textrm{if}\quad \sum_{i=1}^\infty\frac1{\bar b_i\rho_i}=\infty,
\end{array}\right.
\end{equation}
where $\rho_0=0$ and $\rho_i=\Pi_{k=1}^i\frac{\bar d_k}{\bar b_k}$, ($i\in\N$).
\end{thm}

\section{Example and numerical illustrations}\label{bsec:ex}

\begin{figure}
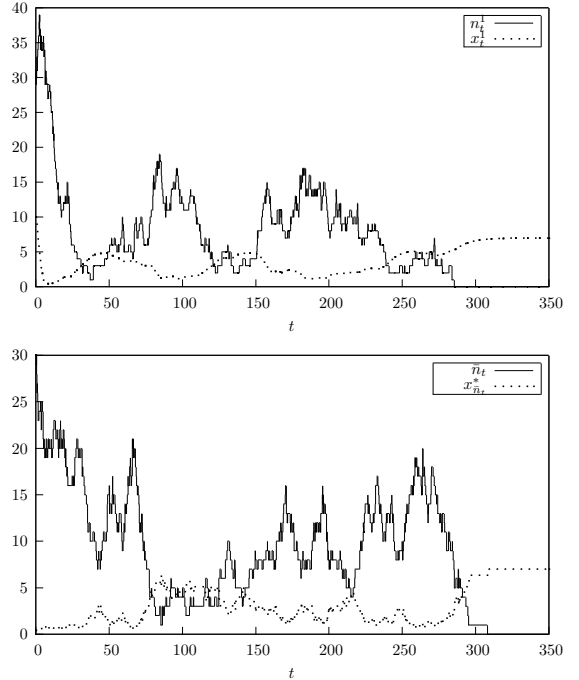

\begin{center}
\scalebox{.6}{\input{tex_traj}}
\scalebox{.6}{\input{tex_traj_av}}
\caption{This figure displays the trajectories of the non-averaged process ${(x^1_t,n^1_t)}_{t\in[0,350]}$ ($\e=1$) starting from $(x^1_0,n^1_0)=(10,30)$ (top) and the averaged process ${(x^*_{\bar n_t},\bar n_t)}_{t\geq0}$ starting from $\bar n_0=30$ (bottom). }\label{fig:traj:m30}
\end{center}
\end{figure}

\begin{figure}
\begin{center}
\includegraphics[width=8cm]{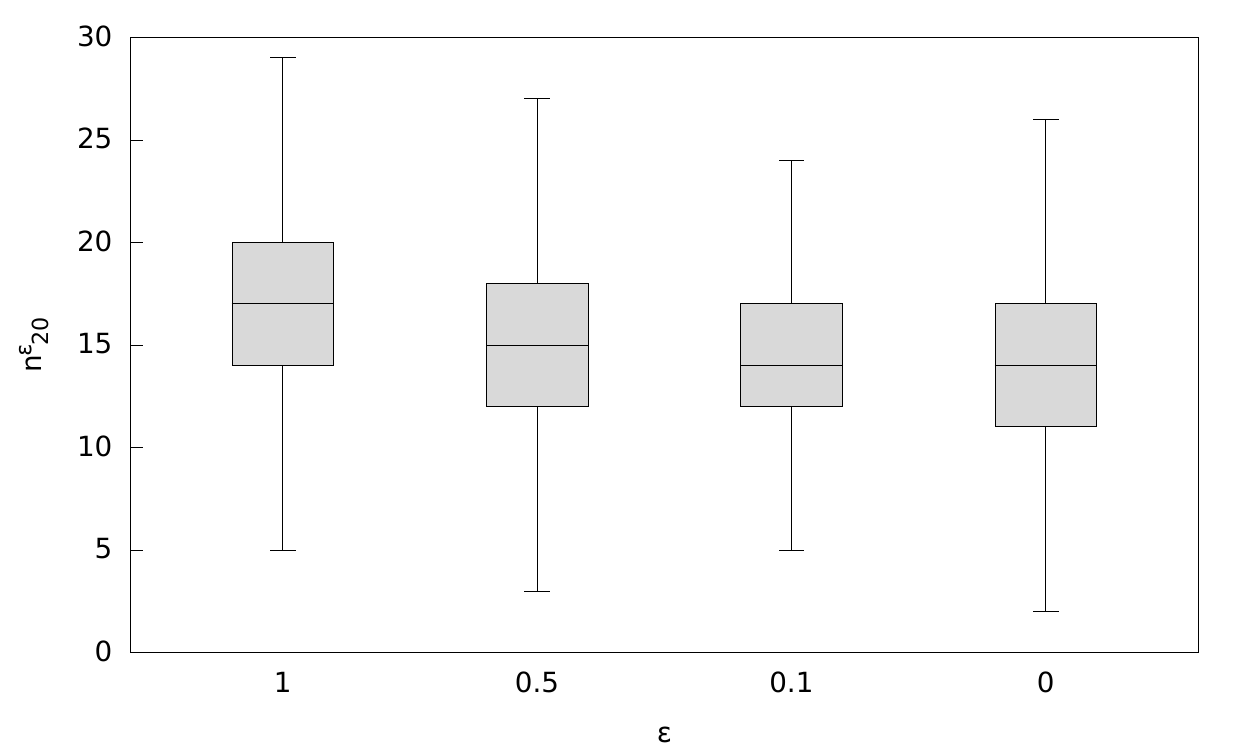}
\\~\\
{\small\begin{tabular}{|c|c|c|c|c|}\hline
$\e$&$1$&$0.5$&$0.1$&$0$\\
\hline
$M_{30}(n^\e_{20})$&$17.0500$&$14.9150$&$14.5740$&$14.2390$\\
\hline
\end{tabular}}
\caption{This figure displays the boxplots of $n^\e_{20}$ starting from $n^\e_0=30$ over $100$ replications for $\e\in\{1,0.5,0.1,0\}$ where $0$ stands here for the averaged model. Bellow the boxplots is a table gathering the empirical mean $M_{30}(n^\e_{20})$ obtained from the same $100$ replications for $\e\in\{1,0.5,0.1,0\}$.}\label{fig:bp:n}
\end{center}
\end{figure}

\begin{figure}
\begin{center}
\includegraphics[width=8cm]{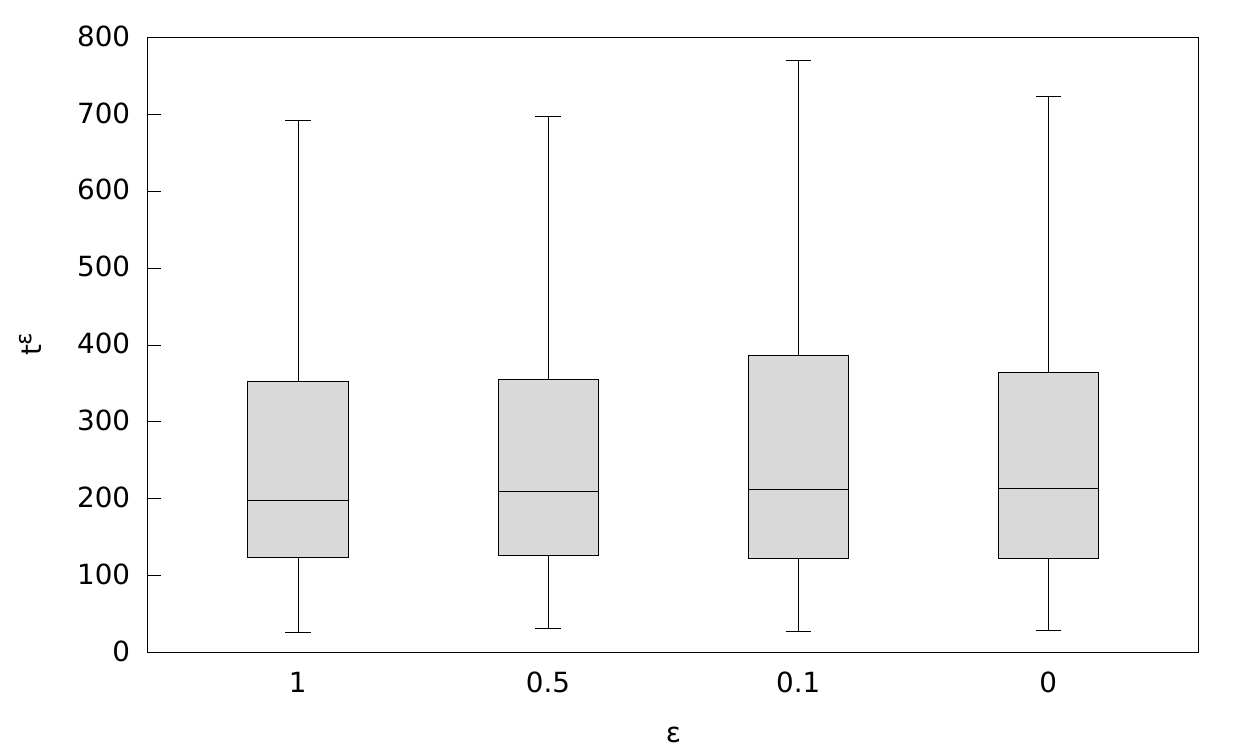}
\\~\\
{\small\begin{tabular}{|c|c|c|c|c|}\hline
$\e$&$1$&$0.5$&$0.1$&$0$\\
\hline
$M_{30}(t^\e)$&$268.1247$&$279.6541$&$287.0660$&$290.2661$\\
\hline
\end{tabular}}
\caption{This figure displays the boxplots of the absorption time $t^\e$ starting from $n^\e_0=30$ over $100$ replications for $\e\in\{1,0.5,0.1,0\}$ where $0$ stands for the averaged model. Bellow the boxplots are a table gathering the empirical mean obtained from the same $100$ replications for $\e\in\{1,0.5,0.1,0\}$.}\label{fig:bp:t}
\end{center}
\end{figure}

In this section, we aim to illustrate the results presented in the previous section. For this purpose, let us consider the model considered in Section \ref{sec:HPP} whose kinetic is given by the equations (\ref{eq:xe}-\ref{eq:ne}) and with time scale parameters $\e\in(0,1)$. As we intend to do some simulation experiments, the constants in the model are fixed as given in Figure \ref{fig:cte}.
\begin{figure}
\begin{center}
\scalebox{.5}{\input{mu}}\\
~\\
\begin{tabular}{|c|c|c|c|c|c|}
\hline
$x_{in}$&$D$&$V$&$\alpha$&$\beta$&$\gamma$\\
\hline
$7$&$0.1$&$1$&$0.5$&$1$&$1$\\
\hline
\end{tabular}
\caption{Numerical values used for the simulations of the model satisfying the equations (\ref{eq:xe2}-\ref{eq:ne2}). At the top is the plot of the function $\mu(x)=\frac{0.15x}{1+x}$.}\label{fig:cte}
\end{center}
\end{figure}\\
For the sake of clarity, let us re-write the system (\ref{eq:xe}-\ref{eq:ne}) with the particular data used in this section,
\begin{equation}\label{eq:xe2}
\frac{\dd x^\e_t}{dt}=0.1(7-x^\e_t)-\frac{0.15}{2}\frac{x^\e_t}{1+x^\e_t}n^\e_t,\quad t\geq0,
\end{equation}
\begin{align}
&\PP(n^\e_{t+h}=n_2|n^\e_t=n_1)\nonumber\\
=&\left\{\begin{array}{ll}
\frac{0.15x^\e_t}{1+x^\e_t}n_1h+\oo(h)&\textrm{if } n_2=n_1+1,n_1\geq1,\\
0.1n_1h+\oo(h)&\textrm{if } n_2=n_1-1,n_1\geq1,\\
1-(\frac{0.15x^\e_t}{1+x^\e_t}+0.1)n_1\hspace{-0.3cm}&h +\oo(h)\\
&\textrm{if } n_2=n_1,n_1\geq1,\\
\oo(h)&\text{else},
\end{array}\label{eq:ne2}\right.
\end{align}
endowed with the initial conditions $x^\e_0$, $n^\e_0$, both deterministic and positives. One can show that, starting from some positive initial value $x^\e_0$, the process ${(x^\e_t)}_{t\geq0}$ stays positive almost-surely. In this case, for any $n\in\N_0$, the quantity $x^*_n$ of Assumption \ref{hyp:1} is the unique positive zero of the function $g(\cdot,n)$ where 
$$
g(x,n)=D(x_{in}-x)-\frac{\alpha}{V}\mu(x)n,\quad x\in\R_+.
$$
Thus, within this particular example, $x^*_n$ is given by
$$
x^*_n=\frac12\left(\sqrt{u^2_n+4v}-u_n\right),
$$ 
where $u_n=7.5n-6$ and $v=7$, for any $n\in\N_0$. Note also that $g(0,n)=Dx_{in}$ for any $n\in\N_0$ such that the sequence ${(g(0,n))}_{n\in\N_0}$ is positive and bounded uniformly in $n\in\N_0$. Moreover, the dissipativity constant $\delta$ is given, with our particular data, by $\delta=D=0.1$. The verification of the other assumptions (Assumption \ref{ass:2} and Assumption \ref{hyp:growth}) are left to the reader.\\
A trajectory of the one time scale process ${(x^1_t,n^1_t)}_{t\geq0}$ ($\e=1$) starting from $(x^1_0,n^1_0)=(10,30)$ is displayed in Figure \ref{fig:traj:m30} (at the top). On this figure, the coupling between $x^1$ and $n^1$ is clearly visible: the growth of the one is correlated to the decrease of the other. We also observe, with the set of data specified in Figure \ref{fig:cte}, the absorption of $n^1$ after some times, corresponding to the prey extinction. Probability and time of extinctions will be discussed later in this section.\\
Applying Theorem \ref{thm:main}, for any finite time horizon $T$, the process ${(n^\e_t)}_{t\in[0,T]}$ converges in law when $\e$ goes to zero towards the birth and death process ${(\bar n_t)}_{t\in[0,T]}$ on $\N_0$ with birth and death rates given by
$$
\bar{b}_n=n\mu(x^*_n)\quad\text{and}\quad \bar{d}_n=Dn,
$$
for $n\in\N$ and $\bar{b}_0=0$ (absorption at zero).\\
A trajectory of the averaged process ${(x^*_{\bar n_t},\bar n_t)}_{t\geq0}$ starting from $\bar n_0=30$ is displayed in Figure \ref{fig:traj:m30} (at the bottom) just below the trajectory of the corresponding non-averaged process. These two processes seem to have qualitatively the same behavior. In particular, in Figure \ref{fig:traj:m30}, both trajectories (averaged or not) are absorbed.
\\
To illustrate the convergence of the slow-fast process ${(x^\e_t,n^\e_t)}_{t\in[0,T]}$ towards the averaged process, we present in Figure \ref{fig:bp:n}, for decreasing values of $\e$, the boxplots of $n^\e_{20}$ and $\bar{n}_{20}$ over $100$ replications when the initial conditions are $n^\e_0=\bar{n}_0=30$ and $x^\e_0=10$. For the same values of $\e$, the empirical means associated to the boxplots, denoted $M_{30}(n^\e_{20})$ and approximating the expectations $\EE_{30}(n^\e_{20})$, are gathered in the table below the boxplots displayed in Figure \ref{fig:bp:n}. We observe, on both the boxplots and the empirical means, the fast convergence of the slow-fast process towards the averaged one when $\e$ goes to zero.
\\
For the set of data used in the present section, one can compute explicitly, using Theorem \ref{thm:abs}, the absorption probability for the averaged process $\bar n$. Indeed, recalling for any $i\in\N$ the expression of $\rho_i$ given in Theorem \ref{thm:abs}, we have
$$
\rho_i=\prod_{k=1}^i\frac{\bar d_k}{\bar b_k}=\frac{D^i}{\prod_{k=1}^i\mu(x^*_i)}.
$$
From the above expression, it is not hard to see that the series $\sum\rho_i$ is divergent such that, according to Theorem \ref{thm:abs}, the averaged process is absorbed with probability one. We can therefore consider its absorption time and compare its mean value to the absorption time for the non-averaged process.\\
In figure \ref{fig:bp:t} are displayed, for decreasing values of $\e$, the boxplots of the absorption time $t^\e$ for the slow-fast process $n^\e$ starting from the initial value $(x^\e_0,n^\e_0)=(10,30)$ over $100$ replications. These boxplots are compared to the boxplot of the absorption time for the averaged process $\bar n$ starting from the same initial value, still over $100$ replications. The empirical means for these different absorption times are also given. At first sight, the displayed boxplots seem comparable whatever the value of $\e$: it could mean that the absorption times for the two-time scales and the averaged models are similar. The convergence of the absorption times is clearer in the table given the empirical mean $M_{30}(t^\e)$ of $t^\e$ starting from the initial value $30$: $M_{30}(t^\e)$ numerically converges towards $M_{30}(t^0)$ when $\e$ goes to zero (the notation $t^0$ standing here for the absorption time of the averaged process).
\\
All the numerical computations have been performed in langage C with the help of the Gnuplot software to obtain graphical outputs. The numerical method that we used to simulate the different trajectories of the considered hybrid models is presented in \cite{R13} as well as convergence results.  

\appendix

\section{Proof of Theorem \ref{thm:main}}\label{sec:pf}

Let $T$ be a finite time horizon. According to (\ref{g:main}) and following \cite[][Ch. 4]{EK09}, for any $\e\in(0,1)$, the process ${(n^\e_t)}_{t\in[0,T]}$ has the following semi-martingale representation
\begin{equation}\label{n:mr}
n^\e_t=n^\e_0+\int_0^t b(x^\e_s,n^\e_s)-d(x^\e_s,n^\e_s)\dd s +M^\e_t,\quad t\geq0,\PP\text{-}a.s.,
\end{equation}
where $M^\e$ is a c\`adl\`ag squared integrable martingale with initial value $M^\e_0=0$. Its quadratic variation is given by
\begin{equation}\label{M:qv}
[M^\e]_t=\int_0^t b(x^\e_s,n^\e_s)+d(x^\e_s,n^\e_s)\dd s,\quad t\geq0.
\end{equation}
The lemma below is about the boundedness of the process ${(n^\e_t)}_{t\in[0,T]}$.
\begin{lem}\label{lem:b}
If $\EE((n^\e_0)^2)$ is bounded uniformly in $\e$, there exists a constant $C_T$, depending on $T$ but otherwise not on $\e$, such that
$$
\EE\left(\sup_{t\in[0,T]}{(n^\e_t)}^2\right)\leq C_T.
$$
\end{lem}
\begin{pf}
From the representation (\ref{n:mr}), $\PP$-a.s., for any $t\in[0,T]$ we have
\begin{align*}
&{(n^\e_t)}^2\\
&\leq 4(n^\e_0)^2+4\left(\int_0^t b(x^\e_s,n^\e_s)-d(x^\e_s,n^\e_s)\dd s \right)^2+2(M^\e_t)^2\\
&\leq 4(n^\e_0)^2+4T\int_0^t \left(b(x^\e_s,n^\e_s)-d(x^\e_s,n^\e_s)\right)^2\dd s +2(M^\e_t)^2\\
&\leq 4(n^\e_0)^2+4T\int_0^t c^2_1+c^2_2(n^\e_s)^2\dd s +2(M^\e_t)^2,
\end{align*}
where we have used successively the fact that the process ${(n^\e_t)}_{t\geq0}$ is $\PP$-a.s. positive, the Jensen inequality, the Assumption \ref{hyp:growth} and the elementary inequality $(a+b)^2\leq 2(a^2+b^2)$. Using $\sup_{u\in[0,s]}(n^\e_u)^2$ as an upper bound for $(n^\e_s)^2$ and taking the supremum over $t\in[0,T]$ in the above inequality we obtain, $\PP$-a.s.,
\begin{align*}
&\sup_{t\in[0,T]}(n^\e_t)^2\\
&\leq 4(n^\e_0)^2+4T\int_0^T c^2_1+c^2_2\sup_{u\in[0,s]}(n^\e_u)^2\dd s +2\sup_{t\in[0,T]}(M^\e_t)^2.
\end{align*}
With the help of the Burkholder-Davis-Gundy inequality: $\EE(\sup_{t\in[0,T]}|M^\e_t|^2)\leq c_4\EE([M^\e]_T)=c_4\int_0^T b(x^\e_s,n^\e_s)+d(x^\e_s,n^\e_s)\dd s$ (where $c_4$ is a positive constant), after some elementary algebra, we have
\begin{align*}
&\EE\left(\sup_{t\in[0,T]}(n^\e_t)^2\right)\\
&\leq 4\EE((n^\e_0)^2)+4T^2c^2_1+2c_4c_1T\\
&+2c_4c_2\int_0^T\sup_{u\in[0,s]}n^\e_u\dd s+4Tc^2_2\int_0^T\sup_{u\in[0,s]}(n^\e_u)^2\dd s.
\end{align*}
A Gronwall-like comparison inequality leads to the existence of a constant $C_T$, only depending on $T$, such that $\EE\left(\sup_{t\in[0,T]}(n^\e_t)^2\right)\leq C_T$.
\end{pf}
\begin{lem}\label{lem:c}
Let $\eta>0$. There exists $\kappa>0$ such that for any $\e>0$ and any stopping times $\sigma,\tau$ such that $\sigma\leq\tau\leq(\sigma+\kappa)\wedge T$, we have
$$
\EE\left((n^\e_\tau-n^\e_\sigma)^2\right)\leq\eta.
$$
\end{lem}
\begin{pf}
Using the optional stopping time theorem and the It\^o isometry, we notice that
\begin{align*}
\EE((M^\e_\tau-M^\e_\sigma)^2)&=\EE((M^\e_\tau)^2)-\EE((M^\e_\sigma)^2)\\
&=\int_\sigma^\tau b(x^\e_s,n^\e_s)+d(x^\e_s,n^\e_s)\dd s.
\end{align*}
Then, the representation (\ref{n:mr}) leads to
\begin{align*}
&\EE((n^\e_\tau-n^\e_\sigma)^2)\\
&=\EE\left(\left(\int_\sigma^\tau b(x^\e_s,n^\e_s)-d(x^\e_s,n^\e_s)\dd s+M^\e_\tau-M^\e_\sigma\right)^2\right)\\
&\leq 2\left(\int_\sigma^\tau b(x^\e_s,n^\e_s)+d(x^\e_s,n^\e_s)\dd s\right)^2\\
&~~~+2\int_\sigma^\tau b(x^\e_s,n^\e_s)+d(x^\e_s,n^\e_s)\dd s.
\end{align*}
The Jensen inequality, Assumption \ref{hyp:growth} and Lemma \ref{lem:b} lead to, after some elementary algebra,
\begin{align*}
&\EE\left((n^\e_\tau-n^\e_\sigma)^2\right)\\
&\leq 2(c^2_1+c^2_2C_T)(\tau-\sigma)^2+2(c_1+c_2C_T)(\tau-\sigma)\\
&\leq C'_T\kappa,
\end{align*}
where $C'_T$ is a constant depending only on $T$. Then, choosing $\kappa$ such that $C'_T\kappa\leq\eta$ yields the result.
\end{pf}
Lemma \ref{lem:b} and \ref{lem:c} leads to the tightness of the family $\{n^\e,\e\in(0,1)\}$ in the Skorohod space $\mathcal{D}([0,T],\N_0)$. We now proceed to the identification of the limit. The proof of the following Lemma is quite similar to the one of Proposition \ref{cv:fast} (recall that $x^\e_0$ is assumed to be positive, $\PP$-a.s.).
\begin{lem}\label{lem:x:b}
The process ${(x^\e_t)}_{t\geq0}$ satisfies
$$
x^\e_t\leq\left(x^\e_0+\frac{1}{\delta}\sup_{n\in\N_0}g(0,n)\right)e^{-\frac{\delta}{\e}t}+\frac{1}{\delta}\sup_{n\in\N_0}g(0,n),
$$
for any $t\geq0$, $\PP$-a.s.
\end{lem}

\begin{lem}\label{lem:last}
Let $\phi$ be a bounded real measurable function on $\N_0$. For any $T>0$, there exists a constant $C_T$ depending only on $T$, such that for any $\e\in(0,1)$ and $t\in[0,T]$,
\begin{itemize}
\item[i)] $\EE\left(\int_0^t[\phi(n^\e_s+1)-\phi(n^\e_s)](b(x^\e_s,n^\e_s)-b(x^*_{n^\e_s},n^\e_s))\dd s\right)$\\$\leq C_T\e$,
\item[ii)] $\EE\left(\int_0^t[\phi(n^\e_s-1)-\phi(n^\e_s)](d(x^\e_s,n^\e_s)-d(x^*_{n^\e_s},n^\e_s))\dd s\right)$\\$\leq C_T\e$.
\end{itemize}
\end{lem}
\begin{pf}
We only prove the statement i) since the proof for the statement ii) is quite similar. Let us denote by ${(T^\e_k)}_{k\geq0}$ the sequence of jumps of the process ${(n^\e_t)}_{t\geq0}$ with $T^\e_0=0$. For any $t\in[0,T]$, splitting the integral according to the sequence ${(T^\e_k)}_{k\geq0}$ and considering all the possible values taken by $n^\e$, we have
\begin{align*}
&\EE\left(\int_0^t[\phi(n^\e_s+1)-\phi(n^\e_s)](b(x^\e_s,n^\e_s)-b(x^*_{n^\e_s},n^\e_s))\dd s\right)\\
&\leq\sum_{k\geq0}\sum_{n\geq0}\EE\bigg(1_{n^\e_{T^\e_k}=n}\int_{T^\e_k\wedge t}^{T^\e_{k+1}\wedge t}|\phi(n+1)-\phi(n)|\bigg.\\
&\hspace{5cm}\bigg.|b(x^\e_s,n)-b(x^*_{n},n)|\dd s\bigg).
\end{align*}
Using Proposition \ref{cv:fast}, Assumption \ref{hyp:growth} and the fact that $\phi$ is bounded, for any $k\geq0$ and $s\in[T^\e_k,T^\e_{k+1})$ we obtain
\begin{align*}
&\EE\left(\int_0^t[\phi(n^\e_s+1)-\phi(n^\e_s)](b(x^\e_s,n^\e_s)-b(x^*_{n^\e_s},n^\e_s))\dd s\right)\\
&\leq2\|\phi\|_\infty\sum_{k\geq0}\sum_{n\geq0}\EE\bigg(1_{n^\e_{T^\e_k}=n}\int_{T^\e_k\wedge t}^{T^\e_{k+1}\wedge t}c_1(1+c_2n+c_3n^2)\bigg.\\
&\hspace{5cm}\bigg.|x^\e_{T^\e_k}-x^*_n|e^{-\frac{\delta}{\e}s}\dd s\bigg).
\end{align*}
Then, Lemma \ref{lem:x:b} and Assumption \ref{hyp:1} yield
\begin{align*}
&\EE\left(\int_0^t[\phi(n^\e_s+1)-\phi(n^\e_s)](b(x^\e_s,n^\e_s)-b(x^*_{n^\e_s},n^\e_s))\dd s\right)\\
&\leq2\|\phi\|_\infty(C_0+3\frac1\delta\sup_{n\in\N_0}|g(0,n)|)\\
&\hspace{1cm}c_1\EE\left(1+c_2\sup_{s\in[0,T]}n^\e_s+c_3\sup_{s\in[0,T]}(n^\e_s)^2\right)\int_0^te^{-\frac{\delta}{\e}s}\dd s,
\end{align*}
where $C_0$ is an almost sure bound for $x^\e_0$. The result follows by Lemma \ref{lem:b} and the fact that $\int_0^te^{-\frac{\delta}{\e}s}\dd s\leq \frac\e\delta$ for any $t\geq0$.
\end{pf}
Since the family $\{n^\e,\e\in(0,1)\}$ is tight, up to the extraction of a subsequence, one may assume that $n^\e$ converges in law when $\e$ goes to zero towards some c\`adl\`ag process $\bar n$ that we intend to characterize. According to the Dynkin formula, for any bounded real measurable function $\phi$ on $\N_0$ and any $t\in[0,T]$, we have
\begin{align*}
&\EE(\phi(n^\e_t))\\
&=\EE(\phi(n^\e_0))+\EE\left(\int_0^t[\phi(n^\e_s+1)-\phi(n^\e_s)]b(x^\e_s,n^\e_s)\dd s\right)\\
&~~~+\EE\left(\int_0^t[\phi(n^\e_s-1)-\phi(n^\e_s)]d(x^\e_s,n^\e_s)\dd s\right).
\end{align*}
Thus, we may consider the decomposition
$$
\EE(\phi(n^\e_t))=I^\e_0+I^\e_b+I^\e_d+J^\e_b+J^\e_d,
$$
where
\begin{align*}
I^\e_0&=\EE(\phi(n^\e_0)),\\
I^\e_b&=\EE\left(\int_0^t[\phi(n^\e_s+1)-\phi(n^\e_s)]b(x^*_{n^\e_s},n^\e_s)\dd s\right),\\
I^\e_d&=\EE\left(\int_0^t[\phi(n^\e_s-1)-\phi(n^\e_s)]d(x^*_{n^\e_s},n^\e_s)\dd s\right),\\
J^\e_b&=\EE\left[\int_0^t[\phi(n^\e_s+1)-\phi(n^\e_s)][b(x^\e_s,n^\e_s)-b(x^*_{n^\e_s},n^\e_s)]\dd s\right],\\
J^\e_d&=\EE\left[\int_0^t[\phi(n^\e_s-1)-\phi(n^\e_s)][d(x^\e_s,n^\e_s)-d(x^*_{n^\e_s},n^\e_s)]\dd s\right].
\end{align*}
Since $n^\e$ converges in law towards $\bar n$, it is not difficult to see that $\EE(\phi(n^\e_t))$, $I^\e_0$, $I^\e_b$ and $I^\e_d$ converges towards $\EE(\phi(\bar n_t))$, $I_0$, $I_b$ and $I_d$ where
\begin{align*}
I_0&=\EE(\phi(\bar n_0)),\\
I_b&=\EE\left(\int_0^t[\phi(\bar n_s+1)-\phi(\bar n_s)]b(x^*_{\bar n_s},\bar n_s)\dd s\right),\\
I_d&=\EE\left(\int_0^t[\phi(\bar n_s-1)-\phi(\bar n_s)]d(x^*_{\bar n_s},\bar n_s)\dd s\right).
\end{align*}
The two terms $J^\e_b$ and $J^\e_d$ go to zero with $\e$ according to Lemma \ref{lem:last}. Finally, the averaged limit process $\bar n$ is characterized by the fact that for any bounded real measurable function $\phi$ on $\N_0$,
\begin{align*}
\EE(\phi(\bar n_t))=&~\EE(\phi(\bar n_0))\\
&+\EE\left(\int_0^t[\phi(\bar n_s+1)-\phi(\bar n_s)]b(x^*_{\bar n_s},\bar n_s)\dd s\right)\\
&+\EE\left(\int_0^t[\phi(\bar n_s-1)-\phi(\bar n_s)]d(x^*_{\bar n_s},\bar n_s)\dd s\right).
\end{align*}
The above equation uniquely characterized the averaged process $\bar n$ as a birth and death process on $\N_0$ with respective parameters
$$
\bar b_n=b(x^*_n,n),\quad\textrm{and}\quad \bar d_n=d(x^*_n,n),
$$
for $n\in\N$ and $b_0=0$ so that $0$ is still an absorbing state. The proof of Theorem \ref{thm:main} is now complete.

\section*{Bibliography}

\bibliographystyle{elsarticle-num} 
\bibliography{mabiblio}

\end{document}

%% file: mu.tex
% GNUPLOT: LaTeX picture
\setlength{\unitlength}{0.240900pt}
\ifx\plotpoint\undefined\newsavebox{\plotpoint}\fi
\begin{picture}(1500,900)(0,0)
\sbox{\plotpoint}{\rule[-0.200pt]{0.400pt}{0.400pt}}%
\put(130,131){\makebox(0,0)[r]{ 0}}
\put(150.0,131.0){\rule[-0.200pt]{4.818pt}{0.400pt}}
\put(130,228){\makebox(0,0)[r]{ 0.02}}
\put(150.0,228.0){\rule[-0.200pt]{4.818pt}{0.400pt}}
\put(130,325){\makebox(0,0)[r]{ 0.04}}
\put(150.0,325.0){\rule[-0.200pt]{4.818pt}{0.400pt}}
\put(130,422){\makebox(0,0)[r]{ 0.06}}
\put(150.0,422.0){\rule[-0.200pt]{4.818pt}{0.400pt}}
\put(130,519){\makebox(0,0)[r]{ 0.08}}
\put(150.0,519.0){\rule[-0.200pt]{4.818pt}{0.400pt}}
\put(130,616){\makebox(0,0)[r]{ 0.1}}
\put(150.0,616.0){\rule[-0.200pt]{4.818pt}{0.400pt}}
\put(130,713){\makebox(0,0)[r]{ 0.12}}
\put(150.0,713.0){\rule[-0.200pt]{4.818pt}{0.400pt}}
\put(130,810){\makebox(0,0)[r]{ 0.14}}
\put(150.0,810.0){\rule[-0.200pt]{4.818pt}{0.400pt}}
\put(150,90){\makebox(0,0){ 0}}
\put(150.0,131.0){\rule[-0.200pt]{0.400pt}{4.818pt}}
\put(384,90){\makebox(0,0){ 2}}
\put(384.0,131.0){\rule[-0.200pt]{0.400pt}{4.818pt}}
\put(619,90){\makebox(0,0){ 4}}
\put(619.0,131.0){\rule[-0.200pt]{0.400pt}{4.818pt}}
\put(853,90){\makebox(0,0){ 6}}
\put(853.0,131.0){\rule[-0.200pt]{0.400pt}{4.818pt}}
\put(1087,90){\makebox(0,0){ 8}}
\put(1087.0,131.0){\rule[-0.200pt]{0.400pt}{4.818pt}}
\put(1322,90){\makebox(0,0){ 10}}
\put(1322.0,131.0){\rule[-0.200pt]{0.400pt}{4.818pt}}
\put(150.0,131.0){\rule[-0.200pt]{0.400pt}{175.375pt}}
\put(150.0,131.0){\rule[-0.200pt]{310.520pt}{0.400pt}}
\put(1439.0,131.0){\rule[-0.200pt]{0.400pt}{175.375pt}}
\put(150.0,859.0){\rule[-0.200pt]{310.520pt}{0.400pt}}
\put(794,29){\makebox(0,0){x}}
\put(1179.0,798.0){\rule[-0.200pt]{0.400pt}{9.877pt}}
\put(1179.0,839.0){\rule[-0.200pt]{57.816pt}{0.400pt}}
\put(1419.0,798.0){\rule[-0.200pt]{0.400pt}{9.877pt}}
\put(1179.0,798.0){\rule[-0.200pt]{57.816pt}{0.400pt}}
\put(1279,819){\makebox(0,0)[r]{$\mu(x)$}}
\put(1299.0,819.0){\rule[-0.200pt]{24.090pt}{0.400pt}}
\put(150,131){\usebox{\plotpoint}}
\multiput(150.58,131.00)(0.493,2.875){23}{\rule{0.119pt}{2.346pt}}
\multiput(149.17,131.00)(13.000,68.130){2}{\rule{0.400pt}{1.173pt}}
\multiput(163.58,204.00)(0.493,2.320){23}{\rule{0.119pt}{1.915pt}}
\multiput(162.17,204.00)(13.000,55.025){2}{\rule{0.400pt}{0.958pt}}
\multiput(176.58,263.00)(0.493,1.964){23}{\rule{0.119pt}{1.638pt}}
\multiput(175.17,263.00)(13.000,46.599){2}{\rule{0.400pt}{0.819pt}}
\multiput(189.58,313.00)(0.493,1.646){23}{\rule{0.119pt}{1.392pt}}
\multiput(188.17,313.00)(13.000,39.110){2}{\rule{0.400pt}{0.696pt}}
\multiput(202.58,355.00)(0.493,1.408){23}{\rule{0.119pt}{1.208pt}}
\multiput(201.17,355.00)(13.000,33.493){2}{\rule{0.400pt}{0.604pt}}
\multiput(215.58,391.00)(0.493,1.210){23}{\rule{0.119pt}{1.054pt}}
\multiput(214.17,391.00)(13.000,28.813){2}{\rule{0.400pt}{0.527pt}}
\multiput(228.58,422.00)(0.493,1.052){23}{\rule{0.119pt}{0.931pt}}
\multiput(227.17,422.00)(13.000,25.068){2}{\rule{0.400pt}{0.465pt}}
\multiput(241.58,449.00)(0.493,0.972){23}{\rule{0.119pt}{0.869pt}}
\multiput(240.17,449.00)(13.000,23.196){2}{\rule{0.400pt}{0.435pt}}
\multiput(254.58,474.00)(0.493,0.814){23}{\rule{0.119pt}{0.746pt}}
\multiput(253.17,474.00)(13.000,19.451){2}{\rule{0.400pt}{0.373pt}}
\multiput(267.58,495.00)(0.493,0.734){23}{\rule{0.119pt}{0.685pt}}
\multiput(266.17,495.00)(13.000,17.579){2}{\rule{0.400pt}{0.342pt}}
\multiput(280.58,514.00)(0.493,0.655){23}{\rule{0.119pt}{0.623pt}}
\multiput(279.17,514.00)(13.000,15.707){2}{\rule{0.400pt}{0.312pt}}
\multiput(293.58,531.00)(0.493,0.616){23}{\rule{0.119pt}{0.592pt}}
\multiput(292.17,531.00)(13.000,14.771){2}{\rule{0.400pt}{0.296pt}}
\multiput(306.58,547.00)(0.493,0.536){23}{\rule{0.119pt}{0.531pt}}
\multiput(305.17,547.00)(13.000,12.898){2}{\rule{0.400pt}{0.265pt}}
\multiput(319.00,561.58)(0.497,0.493){23}{\rule{0.500pt}{0.119pt}}
\multiput(319.00,560.17)(11.962,13.000){2}{\rule{0.250pt}{0.400pt}}
\multiput(332.00,574.58)(0.539,0.492){21}{\rule{0.533pt}{0.119pt}}
\multiput(332.00,573.17)(11.893,12.000){2}{\rule{0.267pt}{0.400pt}}
\multiput(345.00,586.58)(0.590,0.492){19}{\rule{0.573pt}{0.118pt}}
\multiput(345.00,585.17)(11.811,11.000){2}{\rule{0.286pt}{0.400pt}}
\multiput(358.00,597.58)(0.652,0.491){17}{\rule{0.620pt}{0.118pt}}
\multiput(358.00,596.17)(11.713,10.000){2}{\rule{0.310pt}{0.400pt}}
\multiput(371.00,607.59)(0.728,0.489){15}{\rule{0.678pt}{0.118pt}}
\multiput(371.00,606.17)(11.593,9.000){2}{\rule{0.339pt}{0.400pt}}
\multiput(384.00,616.59)(0.728,0.489){15}{\rule{0.678pt}{0.118pt}}
\multiput(384.00,615.17)(11.593,9.000){2}{\rule{0.339pt}{0.400pt}}
\multiput(397.00,625.59)(0.824,0.488){13}{\rule{0.750pt}{0.117pt}}
\multiput(397.00,624.17)(11.443,8.000){2}{\rule{0.375pt}{0.400pt}}
\multiput(410.00,633.59)(0.824,0.488){13}{\rule{0.750pt}{0.117pt}}
\multiput(410.00,632.17)(11.443,8.000){2}{\rule{0.375pt}{0.400pt}}
\multiput(423.00,641.59)(0.950,0.485){11}{\rule{0.843pt}{0.117pt}}
\multiput(423.00,640.17)(11.251,7.000){2}{\rule{0.421pt}{0.400pt}}
\multiput(436.00,648.59)(1.123,0.482){9}{\rule{0.967pt}{0.116pt}}
\multiput(436.00,647.17)(10.994,6.000){2}{\rule{0.483pt}{0.400pt}}
\multiput(449.00,654.59)(1.123,0.482){9}{\rule{0.967pt}{0.116pt}}
\multiput(449.00,653.17)(10.994,6.000){2}{\rule{0.483pt}{0.400pt}}
\multiput(462.00,660.59)(1.214,0.482){9}{\rule{1.033pt}{0.116pt}}
\multiput(462.00,659.17)(11.855,6.000){2}{\rule{0.517pt}{0.400pt}}
\multiput(476.00,666.59)(1.123,0.482){9}{\rule{0.967pt}{0.116pt}}
\multiput(476.00,665.17)(10.994,6.000){2}{\rule{0.483pt}{0.400pt}}
\multiput(489.00,672.59)(1.378,0.477){7}{\rule{1.140pt}{0.115pt}}
\multiput(489.00,671.17)(10.634,5.000){2}{\rule{0.570pt}{0.400pt}}
\multiput(502.00,677.59)(1.378,0.477){7}{\rule{1.140pt}{0.115pt}}
\multiput(502.00,676.17)(10.634,5.000){2}{\rule{0.570pt}{0.400pt}}
\multiput(515.00,682.59)(1.378,0.477){7}{\rule{1.140pt}{0.115pt}}
\multiput(515.00,681.17)(10.634,5.000){2}{\rule{0.570pt}{0.400pt}}
\multiput(528.00,687.60)(1.797,0.468){5}{\rule{1.400pt}{0.113pt}}
\multiput(528.00,686.17)(10.094,4.000){2}{\rule{0.700pt}{0.400pt}}
\multiput(541.00,691.60)(1.797,0.468){5}{\rule{1.400pt}{0.113pt}}
\multiput(541.00,690.17)(10.094,4.000){2}{\rule{0.700pt}{0.400pt}}
\multiput(554.00,695.60)(1.797,0.468){5}{\rule{1.400pt}{0.113pt}}
\multiput(554.00,694.17)(10.094,4.000){2}{\rule{0.700pt}{0.400pt}}
\multiput(567.00,699.60)(1.797,0.468){5}{\rule{1.400pt}{0.113pt}}
\multiput(567.00,698.17)(10.094,4.000){2}{\rule{0.700pt}{0.400pt}}
\multiput(580.00,703.60)(1.797,0.468){5}{\rule{1.400pt}{0.113pt}}
\multiput(580.00,702.17)(10.094,4.000){2}{\rule{0.700pt}{0.400pt}}
\multiput(593.00,707.61)(2.695,0.447){3}{\rule{1.833pt}{0.108pt}}
\multiput(593.00,706.17)(9.195,3.000){2}{\rule{0.917pt}{0.400pt}}
\multiput(606.00,710.61)(2.695,0.447){3}{\rule{1.833pt}{0.108pt}}
\multiput(606.00,709.17)(9.195,3.000){2}{\rule{0.917pt}{0.400pt}}
\multiput(619.00,713.60)(1.797,0.468){5}{\rule{1.400pt}{0.113pt}}
\multiput(619.00,712.17)(10.094,4.000){2}{\rule{0.700pt}{0.400pt}}
\multiput(632.00,717.61)(2.695,0.447){3}{\rule{1.833pt}{0.108pt}}
\multiput(632.00,716.17)(9.195,3.000){2}{\rule{0.917pt}{0.400pt}}
\put(645,720.17){\rule{2.700pt}{0.400pt}}
\multiput(645.00,719.17)(7.396,2.000){2}{\rule{1.350pt}{0.400pt}}
\multiput(658.00,722.61)(2.695,0.447){3}{\rule{1.833pt}{0.108pt}}
\multiput(658.00,721.17)(9.195,3.000){2}{\rule{0.917pt}{0.400pt}}
\multiput(671.00,725.61)(2.695,0.447){3}{\rule{1.833pt}{0.108pt}}
\multiput(671.00,724.17)(9.195,3.000){2}{\rule{0.917pt}{0.400pt}}
\multiput(684.00,728.61)(2.695,0.447){3}{\rule{1.833pt}{0.108pt}}
\multiput(684.00,727.17)(9.195,3.000){2}{\rule{0.917pt}{0.400pt}}
\put(697,731.17){\rule{2.700pt}{0.400pt}}
\multiput(697.00,730.17)(7.396,2.000){2}{\rule{1.350pt}{0.400pt}}
\put(710,733.17){\rule{2.700pt}{0.400pt}}
\multiput(710.00,732.17)(7.396,2.000){2}{\rule{1.350pt}{0.400pt}}
\multiput(723.00,735.61)(2.695,0.447){3}{\rule{1.833pt}{0.108pt}}
\multiput(723.00,734.17)(9.195,3.000){2}{\rule{0.917pt}{0.400pt}}
\put(736,738.17){\rule{2.700pt}{0.400pt}}
\multiput(736.00,737.17)(7.396,2.000){2}{\rule{1.350pt}{0.400pt}}
\put(749,740.17){\rule{2.700pt}{0.400pt}}
\multiput(749.00,739.17)(7.396,2.000){2}{\rule{1.350pt}{0.400pt}}
\put(762,742.17){\rule{2.700pt}{0.400pt}}
\multiput(762.00,741.17)(7.396,2.000){2}{\rule{1.350pt}{0.400pt}}
\put(775,744.17){\rule{2.700pt}{0.400pt}}
\multiput(775.00,743.17)(7.396,2.000){2}{\rule{1.350pt}{0.400pt}}
\put(788,746.17){\rule{2.700pt}{0.400pt}}
\multiput(788.00,745.17)(7.396,2.000){2}{\rule{1.350pt}{0.400pt}}
\put(801,748.17){\rule{2.700pt}{0.400pt}}
\multiput(801.00,747.17)(7.396,2.000){2}{\rule{1.350pt}{0.400pt}}
\put(814,750.17){\rule{2.700pt}{0.400pt}}
\multiput(814.00,749.17)(7.396,2.000){2}{\rule{1.350pt}{0.400pt}}
\put(827,751.67){\rule{3.132pt}{0.400pt}}
\multiput(827.00,751.17)(6.500,1.000){2}{\rule{1.566pt}{0.400pt}}
\put(840,753.17){\rule{2.700pt}{0.400pt}}
\multiput(840.00,752.17)(7.396,2.000){2}{\rule{1.350pt}{0.400pt}}
\put(853,755.17){\rule{2.700pt}{0.400pt}}
\multiput(853.00,754.17)(7.396,2.000){2}{\rule{1.350pt}{0.400pt}}
\put(866,756.67){\rule{3.132pt}{0.400pt}}
\multiput(866.00,756.17)(6.500,1.000){2}{\rule{1.566pt}{0.400pt}}
\put(879,758.17){\rule{2.700pt}{0.400pt}}
\multiput(879.00,757.17)(7.396,2.000){2}{\rule{1.350pt}{0.400pt}}
\put(892,759.67){\rule{3.132pt}{0.400pt}}
\multiput(892.00,759.17)(6.500,1.000){2}{\rule{1.566pt}{0.400pt}}
\put(905,761.17){\rule{2.700pt}{0.400pt}}
\multiput(905.00,760.17)(7.396,2.000){2}{\rule{1.350pt}{0.400pt}}
\put(918,762.67){\rule{3.132pt}{0.400pt}}
\multiput(918.00,762.17)(6.500,1.000){2}{\rule{1.566pt}{0.400pt}}
\put(931,763.67){\rule{3.132pt}{0.400pt}}
\multiput(931.00,763.17)(6.500,1.000){2}{\rule{1.566pt}{0.400pt}}
\put(944,765.17){\rule{2.700pt}{0.400pt}}
\multiput(944.00,764.17)(7.396,2.000){2}{\rule{1.350pt}{0.400pt}}
\put(957,766.67){\rule{3.132pt}{0.400pt}}
\multiput(957.00,766.17)(6.500,1.000){2}{\rule{1.566pt}{0.400pt}}
\put(970,767.67){\rule{3.132pt}{0.400pt}}
\multiput(970.00,767.17)(6.500,1.000){2}{\rule{1.566pt}{0.400pt}}
\put(983,768.67){\rule{3.132pt}{0.400pt}}
\multiput(983.00,768.17)(6.500,1.000){2}{\rule{1.566pt}{0.400pt}}
\put(996,770.17){\rule{2.700pt}{0.400pt}}
\multiput(996.00,769.17)(7.396,2.000){2}{\rule{1.350pt}{0.400pt}}
\put(1009,771.67){\rule{3.132pt}{0.400pt}}
\multiput(1009.00,771.17)(6.500,1.000){2}{\rule{1.566pt}{0.400pt}}
\put(1022,772.67){\rule{3.132pt}{0.400pt}}
\multiput(1022.00,772.17)(6.500,1.000){2}{\rule{1.566pt}{0.400pt}}
\put(1035,773.67){\rule{3.132pt}{0.400pt}}
\multiput(1035.00,773.17)(6.500,1.000){2}{\rule{1.566pt}{0.400pt}}
\put(1048,774.67){\rule{3.132pt}{0.400pt}}
\multiput(1048.00,774.17)(6.500,1.000){2}{\rule{1.566pt}{0.400pt}}
\put(1061,775.67){\rule{3.132pt}{0.400pt}}
\multiput(1061.00,775.17)(6.500,1.000){2}{\rule{1.566pt}{0.400pt}}
\put(1074,776.67){\rule{3.132pt}{0.400pt}}
\multiput(1074.00,776.17)(6.500,1.000){2}{\rule{1.566pt}{0.400pt}}
\put(1087,777.67){\rule{3.132pt}{0.400pt}}
\multiput(1087.00,777.17)(6.500,1.000){2}{\rule{1.566pt}{0.400pt}}
\put(1100,778.67){\rule{3.132pt}{0.400pt}}
\multiput(1100.00,778.17)(6.500,1.000){2}{\rule{1.566pt}{0.400pt}}
\put(1113,779.67){\rule{3.373pt}{0.400pt}}
\multiput(1113.00,779.17)(7.000,1.000){2}{\rule{1.686pt}{0.400pt}}
\put(1127,780.67){\rule{3.132pt}{0.400pt}}
\multiput(1127.00,780.17)(6.500,1.000){2}{\rule{1.566pt}{0.400pt}}
\put(1140,781.67){\rule{3.132pt}{0.400pt}}
\multiput(1140.00,781.17)(6.500,1.000){2}{\rule{1.566pt}{0.400pt}}
\put(1153,782.67){\rule{3.132pt}{0.400pt}}
\multiput(1153.00,782.17)(6.500,1.000){2}{\rule{1.566pt}{0.400pt}}
\put(1166,783.67){\rule{3.132pt}{0.400pt}}
\multiput(1166.00,783.17)(6.500,1.000){2}{\rule{1.566pt}{0.400pt}}
\put(1192,784.67){\rule{3.132pt}{0.400pt}}
\multiput(1192.00,784.17)(6.500,1.000){2}{\rule{1.566pt}{0.400pt}}
\put(1205,785.67){\rule{3.132pt}{0.400pt}}
\multiput(1205.00,785.17)(6.500,1.000){2}{\rule{1.566pt}{0.400pt}}
\put(1218,786.67){\rule{3.132pt}{0.400pt}}
\multiput(1218.00,786.17)(6.500,1.000){2}{\rule{1.566pt}{0.400pt}}
\put(1231,787.67){\rule{3.132pt}{0.400pt}}
\multiput(1231.00,787.17)(6.500,1.000){2}{\rule{1.566pt}{0.400pt}}
\put(1179.0,785.0){\rule[-0.200pt]{3.132pt}{0.400pt}}
\put(1257,788.67){\rule{3.132pt}{0.400pt}}
\multiput(1257.00,788.17)(6.500,1.000){2}{\rule{1.566pt}{0.400pt}}
\put(1270,789.67){\rule{3.132pt}{0.400pt}}
\multiput(1270.00,789.17)(6.500,1.000){2}{\rule{1.566pt}{0.400pt}}
\put(1244.0,789.0){\rule[-0.200pt]{3.132pt}{0.400pt}}
\put(1296,790.67){\rule{3.132pt}{0.400pt}}
\multiput(1296.00,790.17)(6.500,1.000){2}{\rule{1.566pt}{0.400pt}}
\put(1309,791.67){\rule{3.132pt}{0.400pt}}
\multiput(1309.00,791.17)(6.500,1.000){2}{\rule{1.566pt}{0.400pt}}
\put(1283.0,791.0){\rule[-0.200pt]{3.132pt}{0.400pt}}
\put(1335,792.67){\rule{3.132pt}{0.400pt}}
\multiput(1335.00,792.17)(6.500,1.000){2}{\rule{1.566pt}{0.400pt}}
\put(1348,793.67){\rule{3.132pt}{0.400pt}}
\multiput(1348.00,793.17)(6.500,1.000){2}{\rule{1.566pt}{0.400pt}}
\put(1322.0,793.0){\rule[-0.200pt]{3.132pt}{0.400pt}}
\put(1374,794.67){\rule{3.132pt}{0.400pt}}
\multiput(1374.00,794.17)(6.500,1.000){2}{\rule{1.566pt}{0.400pt}}
\put(1387,795.67){\rule{3.132pt}{0.400pt}}
\multiput(1387.00,795.17)(6.500,1.000){2}{\rule{1.566pt}{0.400pt}}
\put(1361.0,795.0){\rule[-0.200pt]{3.132pt}{0.400pt}}
\put(1413,796.67){\rule{3.132pt}{0.400pt}}
\multiput(1413.00,796.17)(6.500,1.000){2}{\rule{1.566pt}{0.400pt}}
\put(1400.0,797.0){\rule[-0.200pt]{3.132pt}{0.400pt}}
\put(1426.0,798.0){\rule[-0.200pt]{3.132pt}{0.400pt}}
\put(150.0,131.0){\rule[-0.200pt]{0.400pt}{175.375pt}}
\put(150.0,131.0){\rule[-0.200pt]{310.520pt}{0.400pt}}
\put(1439.0,131.0){\rule[-0.200pt]{0.400pt}{175.375pt}}
\put(150.0,859.0){\rule[-0.200pt]{310.520pt}{0.400pt}}
\end{picture}